\documentclass[10pt]{article}
\usepackage{amssymb}
\usepackage{amsmath}

\usepackage[hypertex]{hyperref}
\usepackage{url}
\newtheorem{theo}{Theorem}[section]
\newtheorem{prop}[theo]{Proposition}
\newtheorem{lemm}[theo]{Lemma}

\newtheorem{rema}[theo]{Remark}
\newtheorem{Defi}[theo]{Definition}

\newtheorem{conj}[theo]{Conjecture}

\newtheorem{claim}[theo]{Claim}
\newtheorem{fact}[theo]{Fact}
\begin{document}
\author{}
\date{}

\newcommand{\cqfd}
{%
\mbox{}%
\nolinebreak%
\hfill%
\rule{2mm}{2mm}%
\medbreak%
\par%
}
\newfont{\gothic}{eufb10}

\title{Symplectic involutions of $K3$ surfaces act \\
trivially on $CH_0$}

\author{Claire Voisin
\\CNRS, Institut de math\'{e}matiques de Jussieu}

 \maketitle \setcounter{section}{-1}
\section{Introduction}
\setcounter{equation}{0}
For a smooth complex projective variety $X$, Mumford has shown in \cite{mumford}
that the triviality of the Chow group $CH_0(X)$, i.e. $CH_0(X)_{hom}=0$,
implies the vanishing of holomorphic forms of positive degree on $X$.
An immediate generalization is the fact that
a $0$-correspondence $\Gamma\in CH^d(Y\times X)$, with $d={\rm dim}\,X$,
which induces the $0$-map
$\Gamma_*:CH_0(Y)_{hom}\rightarrow CH_0(X)_{hom}$ has the property
that the
maps $\Gamma^*:H^{i,0}(X)\rightarrow H^{i,0}(Y)$  vanish for $i>0$.

Bloch's conjecture is a sort of converse to the above statement, but it needs
the introduction of a certain filtration on $CH_0$ groups of smooth projective varieties. The beginning of this conjectural filtration
is
\begin{eqnarray}\label{formula}F^0CH_0(X)=CH_0(X),\,
F^1CH_0(X)=CH_0(X)_{hom},\\
\nonumber
F^2CH_0(X)={\rm Ker}\,({\rm alb}_X:CH_0(X)_{hom}\rightarrow {\rm Alb}(X)).
\end{eqnarray}
As the filtration is supposed to satisfy $F^{k}CH^0(X)=0$ for $k>{\rm dim}\,X$,
we find that for surfaces, the filtration is fully
determined by
(\ref{formula}).

Bloch's conjecture for correspondences with values in surfaces  is then the following:
\begin{conj} \label{conjbloch}Let $S$ be a smooth projective surface, and let $X$ be a smooth
projective variety, $\Gamma\in CH^2(X\times S)$ be a correspondence such that
the
maps $\Gamma^*:H^{i,0}(S)\rightarrow H^{i,0}(X)$  vanish for $i>0$.
Then
$$\Gamma_*: CH_0(X)_{alb}\rightarrow CH_0(S)$$
vanishes, where $CH_0(X)_{alb}:={\rm Ker}\,
({\rm alb}_X:CH_0(X)_{hom}\rightarrow {\rm Alb}(X))=F^2CH_0(X)$.

\end{conj}
This question can be addressed in particular
to finite group actions on surfaces.
A particular case of the conjecture above is the following:
\begin{conj}\label{conjblochauto} Let $G$ be a finite group acting on a smooth projective complex surface
$S$ with $q=0$. Let $\chi:G\rightarrow \{1,-1\}$ be a character.
Assume that $H^{2,0}(S)^\chi=0$. Then $CH_0(S)_{hom}^\chi=0$.
\end{conj}
Here
$$H^{2,0}(S)^\chi:=\{\omega\in H^{2,0}(S),\,g^*\omega=\chi(g)\omega,\,\forall g\in G\},$$
$$
CH_0(S)_{hom}^\chi:=\{z\in CH_0(S)_{hom},\,g^*z=\chi(g)z,\,\forall g\in G\}.
$$
This is indeed the particular case of the conjecture \ref{conjbloch}
applied to the $0$-correspondence
$$\pi_\chi:=\sum_{g\in G}\chi(g) \Gamma_g\in CH^2(S\times S),$$
where
$\Gamma_g\subset S\times S$ is the graph of $g$.

Conjecture \ref{conjblochauto} is proved in \cite{voisinhodgebloch}
in the situation where $S$ is the zero set of a transverse section
of a $G$-invariant vector bundle on any variety $X$ with trivial Chow groups
(that is $CH^*(X)_{hom}\otimes\mathbb{Q}=0$), under the assumption that $E$ has
many $G$-invariant sections. This generalizes our previous work
in \cite{voisinscuola}, where the case of the Godeaux action of $\mathbb{Z}/5\mathbb{Z}$ on the $CH_0$ group of  invariant quintic surfaces was solved.
This also covers the case (already considered in \cite{voisinscuola})
of the action
of the involution $i$ on $\mathbb{P}^3$ acting with two $-1$ eigenvectors and two
$+1$ eigenvectors on homogeneous coordinates, if we take for $S$ a
quartic surface defined by an $i$-invariant equation and we look at the antiinvariant part
of $CH_0(S)$.

In the paper \cite{huybrechtsJEMS}, Huybrechts proved that
a derived autoequivalence of a $K3$ surface $S$ acting as the identity
on $H^*(S,\mathbb{Z})$ acts as the identity on $CH_0(S)$. The next situation to consider is that of a symplectic finite
order automorphism $g$ of a $K3$ surface
$S$. Thus  $g$ is by definition an automorphism
of $S$ such that $g^*\omega=\omega$, where $\omega$ is the holomorphic $2$-form on $S$.  Such a $g$ acts trivially on $H^{2,0}(S)$ so it has trivial action
on the transcendental lattice of $S$, and the difference
$$g^*-Id\in {\rm Aut}\, H^*(S,\mathbb{Z})$$
is, at least over $\mathbb{Q}$, induced by the cohomology class of a
cycle in $S\times S$ of the form $\sum_i\alpha_iC_i\otimes C'_i$,
where $C_i,\,C'_i$ are curves on $S$ and $\alpha_i$ are rational coefficients.
It seems that if one could take the $\alpha_i$ to be integers,
the above mentioned result of Huybrechts would apply to
show that $g_*$ is the identity on $CH_0(S)$.
Still the problem remains open for these symplectic automorphisms and was explicitly asked by
Huybrechts in \cite{huybrechtssurvey}.
In this note, the case of a symplectic involution $i$ acting on a $K3$ surface
$S$ is considered.
The fact that such symplectic involutions  act trivially on $CH_0(S)$ has been proved
on one hand  in a finite number of  cases in \cite{gulti}, \cite{voisinscuola},
\cite{voisinhodgebloch}, and on the other hand (and more significantly), it has been established in \cite{huybrechtsnew}
for any $K3$ surface with symplectic involution in one of the three series
introduced by van Geemen and Sarti \cite{sartivangeemen} (each series contains itself
an infinite number of families indexed by an integer $d$, and the three series
differ first of all  by the parity of this integer $d$, and secondly, when  $d$ is even, by the  structure of the N\'{e}ron-Severi lattice of the general such surface admitting an invariant line bundle of self-intersection $2d$).

The present paper solves the problem  in general :

\begin{theo}\label{maintheo} Let $S$ be an algebraic $K3$ surface, and let
$i:S\rightarrow S$ be a symplectic involution.
Then $i_*$ acts as the identity on $CH_0(S)$.
\end{theo}
The proof is elementary : It uses the fact that
Prym varieties of \'{e}tale double covers of curves of genus $g$ are of dimension
$g-1$. This departure point is the obvious generalization of the starting
point of Huybrechts and Kemeny's work \cite{huybrechtsnew},
who work with elliptic curves and their \'{e}tale double covers.
This observation  is applied to the \'{e}tale double covers of generic smooth ample
curves $C\subset S/i$ and allows us to prove in section \ref{secproof} that the
group of $i$-antiinvariant $0$-cycles on $S$ is finite dimensional in the Roitman sense (the definition is recalled in section \ref{secfindim}).
One then uses a mild generalization
(Theorem \ref{theoroit} established in section \ref{secfindim}) of a fundamental result due to Roitman
(cf. \cite{roitman1})   in order to conclude that the
group of $i$-antiinvariant $0$-cycles on $S$ is in fact trivial.

\vspace{0.5cm}

{\bf Thanks.} {\it  I thank Daniel  Huybrechts for interesting discussions and comments on this paper.}

\section{Finite dimensionality in the sense of Roitman\label{secfindim}}
Let $X$ be a smooth (connected for simplicity) projective variety over $\mathbb{C}$, and let
$P\subset CH_0(X)$ be a subgroup.

\begin{Defi}\label{defifindim} We will say that $P$ is finite dimensional in the Roitman sense if
there exist a (nonnecessarily connected) smooth projective variety $W$, and a correspondence $\Gamma\subset W\times X$ such that $P$ is contained in the set
$\{\Gamma_*(w),\,w\in W\}$.
\end{Defi}
\begin{rema} {\rm As  $P$ is a subgroup and the cycles $\Gamma_*(w)$ have finitely many possible degrees
(depending on the connected component of $W$ to which $w$ belongs), we conclude
that if $P$ is finite dimensional in the Roitman sense, all elements of $P$ have degree $0$ (so $P\subset CH_0(X)_{hom}$ as  $X$ is connected).
}
\end{rema}
The following result is essentially due to Roitman. (It   is in fact due to Roitman in the case
where $M=X$ and ${\rm Im}\,Z_*=CH_0(X)_{hom}$, see also \cite{voisincime}, lecture 5). The proof we give below is slightly different, as it makes use of Proposition
\ref{simple}, while Roitman uses only elementary arguments. The proof given
 here also has the advantage
that it does not need the torsion freeness of
the group ${\rm Ker}\,({\rm alb}_M:CH_0(M)_{hom}\rightarrow{\rm Alb}\,M))$.

 Let $M$ and $X$ be smooth connected projective varieties with $X$  of dimension $d$. Let $Z\in CH^d(M\times X)$
be a correspondence.
\begin{theo} \label{theoroit} Assume that ${\rm Im}\,(Z_*:CH_0(M)\rightarrow CH_0(X))$
is finite dimensional in the Roitman sense. Then the map
$Z_*: CH_0(M)_{hom}\rightarrow CH_0(X)$ factors through the Albanese
morphism ${\rm alb}_M:CH_0(M)_{hom}\rightarrow {\rm Alb}\,M$ of $M$.
\end{theo}
{\bf Proof.} By definition, there exist a smooth projective variety
$W$ and a  correspondence $\Gamma\subset W\times X$ such that ${\rm Im}\,Z_*$ is contained in the set
$\{\Gamma_*(w),\,w\in W\}$.
Let $C\subset M$ be a curve which is a very general  complete intersection
of sufficiently ample hypersurfaces $H_i\subset M$. Then by the Lefschetz theorem on hyperplane sections, the Jacobian $J(C)$ maps surjectively
to ${\rm Alb}(M)$ and the kernel $K(C)$ is an abelian variety.
 We will prove for completeness the following result:
 \begin{prop}\label{simple} When the $H_i$'s are sufficiciently ample and very general, $K(C)$ is a simple abelian variety.
 \end{prop}
We fix now $C$ as above, satisfying the conclusion of Proposition
\ref{simple} and let $j:C\rightarrow M$ be the inclusion, which induces the  morphism $j_*:J(C)=CH_0(C)_{hom}\rightarrow CH_0(M)$.
We  note that by taking the  $H_i$  sufficiently ample, the dimension
of $K(C)$ can be made arbitrarily large, so we may assume
${\rm dim}\,K>{\rm dim}\, W$.

Let $R\subset K(C)\times  W$ be the following set:
$$R=\{(k,w)\in K(C)\times W,\,Z_*(j_*(k))=\Gamma_*(w)\,\,{\rm in}\,\,CH_0(X)\}.$$
It is known (cf. \cite[10.1.1]{voisinbook})
that $R$ is a countable union of closed irreducible algebraic subsets
$R_i$ of
$K(C)\times W$. As ${\rm Im}\,Z_*$ is contained in
the set
$\{\Gamma_*(w),\,w\in W\}$, the union of the images
of the first projections $p_{\mid R_i}:R_i\rightarrow K(C)$ is equal to $K(C)$.
A Baire category argument then shows that there exists an $i$ such that
$$pr_{1\mid R_i}:R_i\rightarrow K(C)$$
is dominating.
It follows in particular that ${\rm dim}\,R_i\geq {\rm dim}\,K(C)>{\rm dim}\,W$.
The fibers of the second projection
$$pr_{2\mid R_i}:R_i\rightarrow  W$$
are thus positive dimensional.
Let $w\in W$, and $F_w\subset K(C)$ be the fiber over $w$. Then $F_w\subset K(C)$ is positive dimensional,
hence it generates $K(C)$ as a group because $K(C)$ is simple. On the other hand,
by definition of $R$, for any
$f\in F_w$, we have $Z_*(j_*(f))=\Gamma_*(w)$ in $CH_0(X)$, and  $\Gamma_*(w)$ is independent of $f\in F_w$. Hence for any $0$-cycle $z$ of $F_w$,
we have $Z_*(j_*(z))={\rm deg}\,z\,\Gamma_*(w)$ and
 it follows then from the fact that $F_w$ generates $K(C)$ as a group that
 $Z_*\circ j_*$ vanishes identically on $K(C)$.

 In order to conclude that
 $Z_*:CH_0(M)_{hom}\rightarrow CH_0(X)$ factors through
 ${\rm Alb}\,M$, we now observe the following:
 For $k$ large enough, there is a connected subvariety
 $M'$ of $M^k\times M^k$ such that ${\rm Ker}\,{\rm alb}_M$ is generated by cycles
 $z_m=z^+-z^-$ with $z^+_m=\sum_{l\leq k}m_l$,
 $z^-_m=\sum_{k+1\leq l\leq 2k}m_l$, where $m=(m_1,\ldots,m_{2k})\in M'$.
 Furthermore, if the $H_i$'s are taken ample enough, a  very general point $m\in M'$ is supported
on a curve $C$ as above which is very general. Thus
the $0$-cycle $z_m=z^+_m-z^-_m$, being supported on $C$ and annihilated by
${\rm alb}_M$, belongs to $j_*(K(C))$,
and applying the previous reasoning, we conclude that $Z_*(z_m)=0$, for $m$ very general in $M'$.

It remains to prove that it is true for any $m\in M'$.
We can use for this the following easy fact (which is proved by reducing to the case
of curves):
\begin{fact}\label{fact} Let $Y$ be a connected complex projective variety.
Let $U\subset Y$ be the complement of a countable union of proper closed
algebraic subsets. Then any $0$-cycle of $Y$ is rationally equivalent in $Y$
to a
$0$-cycle supported on $U$.
\end{fact}
We apply this observation to $Y=M'$ and to the subset
$U\subset M'$ where we already proved that $Z_*(z_m)=0$ to
conclude that  $m\mapsto Z_*(z_m)$ vanishes identically on
$V'$, hence that $Z_*$ vanishes on ${\rm Ker}\,{\rm alb}_M$.

\cqfd
{\bf Proof of Proposition \ref{simple}.} First of all, we reduce the problem to the
case where $M$ is a surface, by replacing $M$ by a smooth complete intersection
of ample hypersurfaces $T=H_1\cap\ldots\cap H_{m-2}$ and recalling that due to the
Lefschetz theorem on hyperplane sections \cite[2.3.2]{voisinbook}, ${\rm Alb}\,M={\rm Alb}\,T$.
Now we take on $T$ a Lefschetz pencil of very ample hypersurfaces
$T_t,\,t\in\mathbb{P}^1$.
Picard-Lefschetz theory has for consequence (see \cite[3.2.3]{voisinbook}) the irreducibility of the
monodromy action $\rho:\pi_1(\mathbb{P}^1_{reg},t_0)\rightarrow {\rm Aut}\,H^1(T_{t_0},\mathbb{Q})_{van}$, where
$$H^1(T_{t_0},\mathbb{Q})_{van}:={\rm Ker}\,(H^1(T_{t_0},\mathbb{Q})\rightarrow H^3(T,\mathbb{Q}))$$
and $\mathbb{P}^1_{reg}$ is the open set of $\mathbb{P}^1$ parameterizing smooth curves.
As we are working with odd degree cohomology, for which the local monodromies
have infinite order, the same proof as in \cite[3.2.3]{voisinbook} shows as well the irreducibility of the action
of any finite index subgroup $\Gamma\subset \pi_1(\mathbb{P}^1_{reg},t_0)$.

Assume by contradiction that for the general curve $T_t$, the abelian variety
$K(C_t)$ is not simple. Then there is a finite cover
$r:D\rightarrow \mathbb{P}^1$, and a proper sub-abelian
fibration $$\mathcal{A}\subset \mathcal{K}_D,$$
where $\mathcal{K}_D\rightarrow D_{reg}$ is the pull-back to
$D_{reg}:=r^{-1}(\mathbb{P}^1_{reg})$ of the family of abelian varieties
$K(C_t),\,t\in \mathbb{P}^1_{reg}$.
This sub-abelian
fibration (taken up to isogenies) corresponds to
a sub-local system $\mathbb{L}$ of the pull-back to $D_{reg}$ of the local system
on $\mathbb{P}^1_{reg}$ with fiber $H^1(C_t,\mathbb{Q})_{van}$.

The monodromy action on $\mathbb{P}^1_{reg}$ being irreducible on any finite index subgroup of $\pi_1(\mathbb{P}^1_{reg},t_0)$, it is irreducible on the image
$r_*(\pi_1(D_{reg},s_0)),\,r(s_0)=t_0$.
This contradicts the existence of $\mathbb{L}$.

\cqfd
In the next section, we will prove the following:
\begin{prop}\label{findim} Let $S$ be an algebraic $K3$ surface, and let
$i:S\rightarrow S$ be a symplectic involution.
Then the antiinvariant part $CH_0(S)^-=\{z\in CH_0(S),\,i_*(z)=-z\}$
is finite dimensional in the Roitman sense.
\end{prop}
{\bf Proof of Theorem \ref{maintheo}}
We apply Theorem \ref{theoroit} to the case where $X=S$, $M=S$ and
$Z$ is the cycle $\Delta_S-{\rm Graph}(i)$. Here $\Delta_S$ is the diagonal of $S$ and
${\rm Graph}(i)$ is the graph of $i$.
Proposition \ref{findim} says that ${\rm Im}\,Z_*$ is finite dimensional in the Roitman sense and Theorem
\ref{theoroit} then tells us  that $Z_*: CH_0(S)_{hom}\rightarrow CH_0(S)_{hom}$ factors through $ {\rm Alb}\,S=0$. Hence $Z_*$ vanishes on
$CH_0(S)_{hom}$. On the other hand, $Z_*$ is multiplication by
$-2$ on  $CH_0(S)^-\subset CH_0(S)_{hom}$ and we thus
 proved that $CH_0(S)^{-}$ is a $2$-torsion group; as $CH_0(S)$ has no torsion by \cite{roitman2}, we conclude that
$CH_0(S)^-=0$. Thus $Z_*=Id$ on $CH_0(S)$.

\cqfd
\section{Proof of Proposition \ref{findim} \label{secproof}}
We start with the following lemma: Let $M,\,X$ be  smooth projective varieties with
${\rm dim}\,X=d$. Let  $\Gamma\in CH^d(M\times X)$ be a correspondence. Each point
$(m_1,\ldots, m_k)\in M^k$ determines an element $\sum_im_i\in CH_0(M)$. Hence we get a map
$$\Gamma_*:M^k\rightarrow CH_0(X).$$
\begin{lemm} \label{lemmascolaire} Assume there is a point $m\in M$ such that
$\Gamma_*(m)=0$ in $CH_0(X)$ and for some integer $g>0$, one has
$\Gamma_*(M^{g-1})=\Gamma_*(M^{g})$ as subsets of $CH_0(X)$.
Then ${\rm Im} \,\Gamma_*$ is finite dimensional in the Roitman sense.
\end{lemm}
{\bf Proof.} Since $\Gamma_*(M^{g-1})=\Gamma_*(M^{g})$, it is obvious by induction that
$\Gamma_*(M^{g-1})=\Gamma_*M^{k}$ for any $k\geq g-1$.
Any cycle $z\in CH_0(M) $ can be written as $z^+-z^-$, where $z^+ $ and $z^-$ are effective cycles, of degree $k^+,\,k^-$. Up to adding the adequate multiples of $m$ to
$z^+$ and $z^-$, which does not change $\Gamma_*z$, we may assume that $k^+=k^-\geq g$. Thus $\Gamma_*(z)=\Gamma_*(z^+)-\Gamma_*(z^-)$, where $\Gamma_*(z^+)$ and $\Gamma_*(z^-)$ belong to $\Gamma_*(M^{k})=\Gamma_*(M^{g-1})$.
Hence we proved that the
correspondence $\Gamma'\in CH^d( M^{2g-2}\times X)$,
defined as
$$\Gamma'=\sum_{i\leq g-1} (pr_i, p_X)^*\Gamma-\sum_{g\leq i\leq 2g-2} (pr_i, p_X)^*\Gamma$$
satisfies
$${\rm Im}\,\Gamma_*=\Gamma'_*(M^{2g-2}).$$
According to Definition \ref{defifindim}, ${\rm Im}\,\Gamma_*$
is finite dimensional in the Roitman sense.
\cqfd
{\bf Proof of Proposition \ref{findim}.}
Let $S$ be a $K3$ surface endowed with a symplectic involution $i$.
The quotient surface $\Sigma=S/i$ is a singular $K3$ surface.
(By blowing-up its singular points, which correspond to the fixed points of $i$, it
becomes a honest $K3$ surface.)
The canonical bundle of $\Sigma$ (or rather $\Sigma_{reg}$) is trivial.
Let $L\in {\rm Pic }\,\Sigma$ be very ample, and let
$2g-2={\rm deg}\,c_1(L)^2$. By triviality of $K_{\Sigma_{reg}}$, $g$ is the genus of the smooth curves in $|L|$. Furthermore, we have ${\rm dim}\,|L|=g$, due to the exact sequence
$$0\rightarrow \mathbb{C}\rightarrow H^0(\Sigma,L)\rightarrow H^0(C,L_{|C})=H^0(C,K_C)\rightarrow 0,$$
which comes from the similar exact sequence on the desingularization
$\widetilde{\Sigma}$ of $\Sigma$, which has $H^1(\widetilde{\Sigma},\mathcal{O}_{\widetilde{\Sigma}})=0$.

Note also that for a smooth ample curve $C\subset \Sigma$, the inverse image
$\widetilde{C}\subset S$ is smooth, connected, and is an \'{e}tale double  cover
of $C$. (Only the connectedness is to be proved, and this follows from the fact
that otherwise each component $C_1,\,C_2$ of $\widetilde{C}\subset S$ has positive self-intersection and $C_1\cdot C_2=0$ since $\widetilde{C}$ is smooth. This contradicts the Hodge index theorem.)

Let $\Gamma\in CH^2( S\times S)$ be the correspondence
$\Delta_S-{\rm Graph}(i)$.
We prove now the following, where $c_S$ is the effective $0$-cycle of degree $1$ introduced in \cite{beauvillevoisin}:
\begin{claim} We have $\Gamma_*(c_S)=0$ and
$\Gamma_*(S^g)=\Gamma_*(S^{g-1})$.
\end{claim}
According to Lemma \ref{lemmascolaire}, this proves Proposition \ref{findim}, since
$CH_0(S)^-={\rm Im}\,\Gamma_*$. (The last fact follows from the fact that $\Gamma_*$
acts as $-2 \,Id$ on $CH_0(S)^-$, which is a divisible group.)

\vspace{0.5cm}

{\bf Proof of the claim.} The cycle $c_S$ is obviously $i$-invariant since it is the class of any point of $S$ belonging to a rational curve $D\subset S$, and
if $x\in D$ then $i(x)\in i(D)$ also belongs to a rational curve in $S$.

Let $s=(s_1,\ldots, s_g)$ be a general point of $S^g$.
Then if we denote by $\sigma_i$ the image of $s_i$ in $\Sigma=S/i$,
the $g$-uple $(\sigma_1,\ldots,\sigma_g)$ is generic in $\Sigma^g$
and there exists a unique curve $C_s\in |L|$ containing all the $\sigma_i$'s.
The curve $C_s$ being general in $|L|$, it is smooth and
thus we have the \'{e}tale double cover $\widetilde{C_s}\rightarrow C_s$, with
$\widetilde{C_s}\subset S$ containing the points $s_i$.
Consider the $0$-cycle
$$z_s=\sum_ls_l-i(\sum_ls_l)=\Gamma_*(\sum_ls_l)\in CH_0(S).$$
This cycle clearly depends only on the Abel image
$${\rm alb}_{\widetilde{C_s}}(\sum_ls_l-i(\sum_ls_l)),$$
which is an antiinvariant element of $J(\widetilde{C_s})$ or, up to
$2$-torsion, an element of the Prym variety $P(\widetilde{C_s}/C_s)$ which is
a $g-1$-dimensional abelian variety.

In other words, we find that, on a Zariski open set $U$ of $S^g$,
the map $$S^g\rightarrow  CH_0(S)^-,\,(s_1,\ldots,s_g)\mapsto z_s,$$
factors through the morphism
$$ f:U\rightarrow \mathcal{P}(\widetilde{\mathcal{C}}/\mathcal{C}),
(s_1,\ldots,s_g)\mapsto {\rm alb}_{\widetilde{C_s}}(s_1+\ldots+s_g-i(s_1)-\ldots - i(s_g)),$$
where $\mathcal{C}\rightarrow |L|_{0}$ is the universal smooth curve over the Zariski open set $|L|_0$ of $|L|$ parameterizing
smooth curves, $\widetilde{\mathcal{C}}\rightarrow |L|_{0}$ is the universal family of double covers, and
$\mathcal{P}(\widetilde{\mathcal{C}}/\mathcal{C})\rightarrow |L|_{0}$
is the corresponding
Prym fibration.

The total space of the Prym fibration $\mathcal{P}(\widetilde{\mathcal{C}}/\mathcal{C})$
has dimension $2g-1$, while $U$ has dimension $2g$, so the morphism
$f$ has positive dimensional fibers. It follows that for
$s\in U$, there is a curve $F_s\subset S^g$
such that the $0$-cycle
$z_t=\sum_lt_l-i(\sum_lt_l)$ is rationally equivalent  to $z_s$ in $S$ for any $(t_1,\ldots, t_l)\in F_s$.
Choose an ample   curve
$D\subset S$ whose irreducible components are rational. The curve
$F_s$ meets the ample divisor
$\sum_lpr_l^{-1}(D)$, where $pr_l:S^g\rightarrow S$ is the $l$-th projection.
Hence the $0$-cycle $z_s$ is rationally equivalent to
a $0$-cycle of the form
$z_t=\sum_lt_l-i(\sum_lt_l)$, where we have $t_{l_0}\in D$ for some $l_0$.
We have seen already that the $0$-cycle $t_{l_0}-i(t_{l_0})$
vanishes in $CH_0(S)$ and it follows that
$z_s$ is rationally equivalent to the cycle $\sum_{l\not=l_0}t_l-i(\sum_lt_{l\not=l_0})$.
Thus $z_s\in \Gamma_*(S^{g-1})$ for $s=(s_1,\ldots , s_g)\in U$.

To conclude the proof, we have to show that the above result is true for any
$(s_1,\ldots , s_g)\in S^g$. This follows from the statement
in Fact \ref{fact}, which
we apply  to $Y=S^g$  to  conclude that the cycles $z_s$ for $s=(s_1,\ldots,s_g)\in U$ fill-in the image $\Gamma_*(S^g)$.
 Proposition \ref{findim} is thus proved.

\cqfd

Institut de Math\'{e}matiques de Jussieu

Equipe Topologie et  G\'{e}om\'{e}trie  alg\'{e}briques

Case 247, 4 Place Jussieu,

75005 Paris, France

voisin@math.jussieu.fr
\end{document}